\documentclass{svjour3}                     
\smartqed  
\usepackage{graphicx}
\usepackage{latexsym}
\usepackage{amssymb}
\usepackage{amsmath}
\usepackage{algorithm} 
\usepackage{algorithmic}
\usepackage[misc]{ifsym}

\newtheorem{them}{Theorem}[section]

\newtheorem{pro}{Proposition}[section]

\journalname{JOGO}
\begin{document}

\title{Convex mixed integer nonlinear programming problems and an outer approximation algorithm}


\titlerunning{Convex MINLP and Outer Approximation Algorithm}



\author{Zhou Wei \and M. Montaz Ali}


\institute{Zhou Wei(\Letter) \at Department of Mathematics, Yunnan University,  Kunming 650091, People's Republic of China\\ \email{wzhou@ynu.edu.cn} \and M. Montaz Ali \at School of Computational and Applied Mathematics, University of the Witwatersrand, \\ Wits 2050, Johannesburg, South Africa, \\ and \\
TCSE, Faculty of Engineering and Built Environment, University of the Witwatersrand,\\ Wits 2050, Johannesburg, South Africa.\\ \email{Montaz.Ali@wits.ac.za}}

\date{Received: date / Accepted: date}

\maketitle

\begin{abstract}
In this paper, we mainly study one class of convex mixed-integer nonlinear programming problems (MINLPs) with non-differentiable data. By dropping the differentiability assumption, we substitute gradients with subgradients obtained from KKT conditions, and use the outer approximation method to reformulate convex MINLP as one equivalent MILP master program. By solving a finite sequence of subproblems and relaxed MILP problems, we establish an outer approximation algorithm to find the optimal solution of this convex MINLP. The convergence of this algorithm is also presented. The work of this paper generalizes and extends the outer approximation method in the sense of dealing with convex MINLPs from differentiable case to non-differentiable one.


\end{abstract}

\section{Introduction}

Many practical optimization problems are modelled as mixed-integer nonlinear programming problems (MINLPs) involving continuous and discrete variables and the study of solution algorithms for these optimization problems has been an active focus of research over the past decades (cf. \cite{B,F,G,GS1,GS2,L1,L2,TS,TS1,W2} and references therein). Suppose $f,g_i: \mathbb{R}^n\times \mathbb{R}^p\rightarrow \mathbb{R}$ $ (i=1,\cdots,m)$ are nonlinear functions, $X$ is a nonempty compact convex set in $\mathbb{R}^n$ and $Y$ is a set of discrete variables in $\mathbb{R}^p$. The general form for MINLPs is defined mathematically as follows:
\begin{equation}\label{1.1}
{\rm (P)} \left \{
\begin{array}l
\mathop{\rm minimize}\limits_{x,\,y} \ \ \ f(x, y)\\
 {\rm subject\ to}\ \  g_i(x, y)\leq 0, i=1,\cdots,m,\\
\ \ \ \ \ \ \ \ \ \ \ \ \ \ \   x\in X, y\in Y\ {\rm discrete \ variable}.
\end{array}
\right.
\end{equation}
This paper is devoted to one class of convex MINLPs in which objective and constraint functions $f,g_i$ for $i=1,\cdots,m$ are convex but not differentiable.

The class of convex MINLPs has been extensively studied by many authors and several methods for these MINLPs have been developed over past decades. These methods include branch-and-bound, generalized Benders decomposition, extended cutting-plane method, NLP/LP based branch and bound and outer approximation method (cf. \cite{B1,DG,EMW,E2,FL,FR,G1,G,L1,NV,WA,W1,W2} and references therein). Note that the extended cutting-plane method was proposed by Westerlund and Pettersson \cite{W1} for solving differentiable convex MINLPs. Subsequently, Westerlund and Pettersson \cite{W2} presented this method to deal with a more general case of MINLPs including pseudo-convex functions. It was shown in \cite{W2} that one MINLP with pseudo-convex functions and pseudo-convex constraints can be solved to global optimality by the cutting-plane techniques. In 2014, Eronen, M\"{a}kel\"{a} and Westerlund \cite{EMW} generalized the extended cutting-plane method for solving convex nonsmooth MINLPs and provided one ECP algorithm which was proved to converge to one global optimum. Recently they \cite{E2} further considered this extended cutting plane method to deal with nonsmooth MINLPs with pseudo-convexity assumptions.

It is known that Duran and Grossmann \cite{DG} introduced the outer approximation method to deal with a particular class of MINLPs which was restricted to contain separable convex differentiable functions and not general convex differentiable functions in all variables. These separable convex functions were composed of convex differentiable functions in continuous variables and linear functions in discrete variables separately. Afterwards Fletcher and Leyffer \cite{FL} further extended the outer approximation method for solving convex MINLPs with convex and continuously differentiable objective constraint functions, and provided a linear outer approximation algorithm to attain the optimal solution of this MINLP by solving a finite sequence of relaxed subproblems. This extension is the pioneering work on outer approximation method in a sense of solving MINLPs where the discrete variables are considered as nonlinear. In 2008, Bonami et. al \cite{B1} also studied outer approximation algorithms for convex and continuously differentiable MINLPs. Recently the authors in \cite{EMW} and \cite{WA} used the outer approximation method to study convex nonsmooth MINLPs and established the resulting algorithms. It is noted that differentiability of functions plays an important role in the construction of relaxation and is proved to be an important matter for allowing to solve these relaxed subproblems efficiently. Since nonsmooth optimization problems defined by non-differentiable functions appear in practice, from the theoretical viewpoint as well as for applications, it is interesting and significant to consider convex and non-differentiable MINLPs. Motivated by this, in this paper, we are inspired by \cite{B1,DG,EMW,FL,WA} to continue studying one convex MINLP by dropping the differentiability assumption and aim to construct an outer approximation algorithm for solving this MINLP. The outer approximation method used herein is along the line given in \cite{FL,WA} and consists of the use of KKT conditions to linearize the objective and constraint functions at different points for constructing an equivalent MILP relaxation of the problem.

The paper is organized as follows. In section 2, we give some definitions and preliminaries used in this paper.  Section 3 contains the equivalent reformulation of convex MINLP by the outer approximation method and one outer approximation algorithm for finding optimal solutions of this MINLP. The reformulation are mainly dependent on KKT conditions and projection techniques. For the algorithm construction, it is necessary to solve a finite sequence of nonlinear programs including feasible and infeasible subproblems and the relaxations of mixed-integer linear master program. The convergence theorem for the established algorithm is also presented therein. The conclusion of this paper is presented in section 4. Section 5 is an Appendix which contains the proofs of the main results given for constructing the algorithm in the paper.

\setcounter{equation}{0}
\section{Preliminaries}
Let $\|\cdot\|$ be the norm of $\mathbb{R}^n$ and denote the inner product between two elements of $\mathbb{R}^n$ by $\langle\cdot, \cdot\rangle$.
Let $\Omega$ be a closed convex set of $\mathbb{R}^n$ and $x\in \Omega$. We denote $T(\Omega, x)$ the  contingent cone of $\Omega$ at $x$; that is, $v\in T(\Omega, x)$ if and only if there exist a sequence $\{v_k\}$ in $\mathbb{R}^n$ converging to $v$ and a sequence $t_k$ in $(0, +\infty)$ decreasing to $0$ such that $x+t_kv_k\in \Omega$ for all $k\in\mathbb{N}$, where $\mathbb{N}$ denotes the set of all natural numbers. It is known from \cite{AF} that
$$
T(\Omega, x)=cl(\mathbb{R}_+(\Omega-x))
$$
where $cl$ denotes the closure.

Let $N(\Omega, x)$ denote the normal cone of $\Omega$ at $x$, that is
\begin{equation}
N(\Omega, x):=\{\gamma\in \mathbb{R}^n: \langle\gamma, z-x\rangle\leq 0\ \ {\rm for\ all} \ z\in \Omega\}.
\end{equation}
It is easy to verify that normal cone $N(\Omega, x)$ and contingent cone $T(\Omega, x)$ are the polar dual; that is
$$
N(\Omega, x)=\big(T(\Omega, x)\big)^{\circ}:=\big\{\gamma\in\mathbb{R}^n: \langle\gamma, v\rangle\leq 0 \ \ {\rm for\ all} \ v\in T(\Omega, x)\big\}.
$$
Readers are invited to consult the book \cite{AF} for more details on contingent cone and normal cone.


Let $\varphi: \mathbb{R}^n\rightarrow\mathbb{R}$ be a continuous convex function, $\bar x\in\mathbb{R}^n$ and $h\in\mathbb{R}^n$. Recall (cf. \cite{P}) that $d^+\varphi(\bar x)(h)$ denotes the directional derivative of $\varphi$
at $\bar x$ along the direction $h$ and is defined by
$$
d^+\varphi(\bar x)(h):=\lim\limits_{t\rightarrow 0^+}\frac{\varphi(\bar x+th)-\varphi(\bar x)}{t}.
$$
We denote $\partial\varphi(\bar x)$ the subdifferential of $\varphi$ at $\bar x$ which is defined by
$$
\partial\varphi(\bar x):=\{\alpha\in \mathbb{R}^n:\; \langle\alpha, x-\bar x\rangle \leq \varphi(x)-\varphi(\bar x)\ {\rm for\ all} \ x\in\mathbb{R}^n\}.
$$
Each vector in $\partial\varphi(\bar x)$ is called a subgradient of $\varphi$ at $\bar x$. It is known from \cite{P} that $\alpha\in\partial\varphi(\bar x)$ if and only if
$$
\langle\alpha, h\rangle\leq d^+\varphi(\bar x)(h)\ \ {\rm for\ all} \ h\in \mathbb{R}^n.
$$

Recall that $\varphi$ is said to be G\^ateaux differentiable at $\bar x$ if there exists $d\varphi(\bar x)\in\mathbb{R}^n$ such that
\begin{equation}\label{2.2a}
\lim\limits_{t\rightarrow 0^+}\frac{\varphi(\bar x+th)-\varphi(\bar x)}{t}=\langle d\varphi(\bar x), h\rangle\ \ {\rm for\ all}\ h\in\mathbb{R}^n
\end{equation}
and $\varphi$ is said to be Fr\'echet differentiable at $\bar x$ if $\varphi$ is G\^ateaux differentiable there and the limit in \eqref{2.2a} exists uniformly for $\|h\|\leq 1$ as $t\rightarrow 0^+$.

It is known from \cite{P} that $\varphi$ is G\^ateaux differentiable at $\bar x$ if and only if $\partial\varphi(\bar x)$ is the singleton. Further, G\^ateaux differentiability of $\varphi$ is equivalent to the Fr\'echet differentiability of $\varphi$ due to the local Lipschitzian property of $\varphi$ and the compactness of unit closed ball in $\mathbb{R}^n$.

Given a continuous convex function $\phi: \mathbb{R}^n\times\mathbb{R}^p\rightarrow \mathbb{R}$ and $(\bar x, \bar y)\in\mathbb{R}^n\times\mathbb{R}^p$, one vector $(\alpha, \beta)\in \mathbb{R}^n\times\mathbb{R}^p$ is the subgradient of $\phi$ at $(\bar x, \bar y)$ if and only if
\begin{equation}\label{2.2}
\phi(x, y)\geq\phi(\bar x, \bar y)+(\alpha, \beta)^T\begin{pmatrix}x-\bar x\\y-\bar y\end{pmatrix}\ \ {\rm for\ all} \ (x, y)\in\mathbb{R}^n\times\mathbb{R}^p,
\end{equation}
where $(\alpha, \beta)^T$ is the transpose of matrix $(\alpha, \beta)$. When $\bar y$ is fixed ({\it resp.} $\bar x$ is fixed), the subdifferential of $\phi(\cdot,\bar y)$ ({\it resp.} $\phi(\bar x,\cdot)$)  at $\bar x$ ({\it resp.} $\bar y$) is the set defined by
\begin{equation*}
\partial \phi(\cdot,\bar y)(\bar x):=\big\{\alpha\in \mathbb{R}^n: \phi(x, \bar y)\geq\phi(\bar x, \bar y)+\langle\alpha, x-\bar x\rangle\ \ {\rm for\ all} \ x\in\mathbb{R}^n\big\}
\end{equation*}
$$
\Big({\it resp.}\  \partial \phi(\bar x,\cdot)(\bar y):=\big\{\beta\in \mathbb{R}^p: \phi(\bar x, y)\geq\phi(\bar x, \bar y)+\langle\beta, y-\bar y\rangle\ \ {\rm for\ all} \ y\in\mathbb{R}^p\big\}\Big).
$$

The following proposition on the subdifferential of convex functions is easy to verify from the definition.

\begin{pro}
Let $\phi: \mathbb{R}^n\times\mathbb{R}^p\rightarrow \mathbb{R}$ be a continuous convex function and $(\bar x, \bar y)\in\mathbb{R}^n\times\mathbb{R}^p$. Then for any $(\alpha, \beta)\in \partial\phi(\bar x, \bar y)$, one has $\alpha\in\partial\phi(\cdot,\bar y)(\bar x)$ and $\beta\in\partial\phi(\bar x,\cdot)(\bar y)$.
\end{pro}

It is an interesting question to consider the converse of Proposition 2.1. This question is also interesting even for smooth convex functions in mathematical analysis. The question is explicitly stated as follows:

{\it Given one continuous convex function $\phi: \mathbb{R}^n\times\mathbb{R}^p\rightarrow \mathbb{R}$ and one vector $\bar\alpha$ from $\partial\phi(\cdot,\bar y)(\bar x)$, whether or not is there some vector $\bar\beta\in \mathbb{R}^p$ such that $(\bar\alpha, \bar\beta)\in \partial\phi(\bar x, \bar y)$?}

The following propositions provided an affirmative answer to this question. These propositions will play a key role in construction of outer approximation algorithm in the sequel. The first proposition is on convex and Fr\'echet differentiable functions.
\begin{pro}
Let $\phi:\mathbb{R}^n\times\mathbb{R}^p\rightarrow \mathbb{R}$ be a continuous convex function and $(\bar x, \bar y)\in \mathbb{R}^n\times\mathbb{R}^p$. Suppose that $\phi(\cdot,\bar y)$ is Fr\'echet differentiable at $\bar x$ and $\phi(\bar x, \cdot)$ is Fr\'echet differentiable at $\bar y$. Then $\phi$ is Fr\'echet differentiable at $(\bar x, \bar y)$.
\end{pro}

\noindent{\bf Proof.} By the Fr\'echet differentiability of $\phi(\cdot,\bar y)$ and $\phi(\bar x, \cdot)$, one has
$$
\partial\phi(\cdot,\bar y)(\bar x)=\{\triangledown_x\phi(\bar x,\bar y)\}\ \ {\rm and} \ \ \partial\phi(\bar x,\cdot)(\bar y)=\{\triangledown_y\phi(\bar x,\bar y)\}.
$$
This and Proposition 2.1 imply that $\partial\phi(\bar x, \bar y)$ is the singleton and
$$
\partial\phi(\bar x, \bar y)=\{(\triangledown_x\phi(\bar x,\bar y), \triangledown_y\phi(\bar x,\bar y)\}.
$$
Hence $\phi$ is G\^{a}teaux differentiable at $(\bar x, \bar y)$ and consequently Fr\'echet differentiable at $(\bar x, \bar y)$. The proof is complete. $\Box$

Proposition 2.2 may not necessarily be true for non-convex functions. Consider function $\phi$ on $\mathbb{R}\times\mathbb{R}$ defined as: $\phi(x,y)=\frac{x^2y^2}{(x^2+y^2)^{3/2}}$ if $x^2+y^2\not=0$ and $\phi(x,y)=0$ if $x^2+y^2=0$. Then $\phi$ is continuous on  $\mathbb{R}\times\mathbb{R}$ and partial derivatives $\triangledown_x\phi(0,0)$ and $\triangledown_y\phi(0,0)$ exist ($\triangledown_x\phi(0,0)=\triangledown_y\phi(0,0)=0$). However, one can verify that $\phi$ is not differentiable at $(0,0)$.

\begin{pro}
Let $\phi:\mathbb{R}^n\times\mathbb{R}^p\rightarrow \mathbb{R}$ be a continuous convex function and $(\bar x, \bar y)\in \mathbb{R}^n\times\mathbb{R}^p$. Then for any $\bar\alpha\in\partial\phi(\cdot,\bar y)(\bar x)$, there exists $\bar\beta\in\mathbb{R}^p$ such that $(\bar\alpha, \bar\beta)\in\partial\phi(\bar x, \bar y)$.
\end{pro}

\noindent{\bf Proof.} Let $F_{\bar y}:\mathbb{R}^n\rightarrow\mathbb{R}^n\times\mathbb{R}^p$ be defined by $F_{\bar y}(x):=(x, \bar y)$. Then $\phi(\cdot,\bar y)=\phi\circ F_{\bar y}$, and it is easy to verify that $F_{\bar y}$ is differentiable at $\bar x$ and
\begin{equation}\label{2.3}
  \triangledown F_{\bar y}(\bar x)(h)=(h, 0)\in \mathbb{R}^n\times\mathbb{R}^p
\end{equation}
holds for all $h\in \mathbb{R}^n$. Let $\bar\alpha\in\partial\phi(\cdot,\bar y)(\bar x)$. We first prove that
\begin{equation}\label{2.4}
  \bar\alpha\in \triangledown F_{\bar y}(\bar x)^*(\partial\phi(\bar x, \bar y))
\end{equation}
where $\triangledown F_{\bar y}(\bar x)^*$ is the conjugate operator of $\triangledown F_{\bar y}(\bar x)$.

Since $\phi$ is continuous at $(\bar x, \bar y)$, it follows that $\partial\phi(\bar x, \bar y)$ is a nonempty, convex and compact subset by \cite[Proposition 1.11]{P} and then $\triangledown F_{\bar y}(\bar x)^*(\partial\phi(\bar x, \bar y))$ is convex and compact as $\triangledown F_{\bar y}(\bar x)^*$ is continuous.

Suppose to the contrary that $\bar\alpha\not\in \triangledown F_{\bar y}(\bar x)^*(\partial\phi(\bar x, \bar y))$. By the seperation theorem, there exists $\bar u\in\mathbb{R}^n$ with $\|\bar u\|=1$ such that
\begin{eqnarray*}
\langle\bar\alpha, \bar u\rangle&>&\max\{\langle\triangledown F_{\bar y}(\bar x)^*(\alpha, \beta), \bar u\rangle: (\alpha, \beta)\in\partial\phi(\bar x, \bar y)\}\\
&=& \max\{\langle(\alpha, \beta), \triangledown F_{\bar y}(\bar x)(\bar u)\rangle: (\alpha, \beta)\in\partial\phi(\bar x, \bar y)\}.
\end{eqnarray*}
This and \eqref{2.3} imply that
\begin{equation}\label{2.5}
\langle\bar\alpha, \bar u\rangle>\max\{\langle(\alpha, \beta), (\bar u, 0)\rangle: (\alpha, \beta)\in\partial\phi(\bar x, \bar y)\}.
\end{equation}
Noting that $\bar\alpha\in\partial\phi(\cdot,\bar y)(\bar x)$ and $\phi$ is a continuous convex function on $\mathbb{R}^n\times\mathbb{R}^p $, it follows from \cite[Proposition 2.24]{P} and \eqref{2.5} that
\begin{equation*}
  d^+\phi(\cdot,\bar y)(\bar x)(\bar u)\geq\langle\bar\alpha, \bar u\rangle>d^+\phi(\bar x, \bar y)(\bar u, 0)=d^+\phi(\cdot,\bar y)(\bar x)(\bar u),
\end{equation*}
which is contradiction. Thus \eqref{2.4} holds.

By virtue of \eqref{2.4}, there exists $(\hat{\alpha}, \bar{\beta})\in\partial\phi(\bar x,\bar y)$ such that $\bar\alpha=\triangledown F_{\bar y}(\bar x)^*(\hat{\alpha}, \bar{\beta})$. It suffices to prove that $\bar\alpha=\hat{\alpha}$.

For any $h\in\mathbb{R}^n$, by using \eqref{2.3}, one has
$$
\langle\bar\alpha, h\rangle=\langle\triangledown F_{\bar y}(\bar x)^*(\hat\alpha, \bar\beta), h\rangle=\langle(\hat\alpha, \bar\beta), \triangledown F_{\bar y}(\bar x)(h)\rangle=\langle(\hat\alpha, \bar\beta), (h, 0)\rangle=\langle\hat\alpha, h\rangle.
$$
This means that $\bar\alpha=\hat\alpha$. The proof is complete. $\Box$

The following proposition is on the subdifferential of maximum function of two convex functions which is from \cite[Theorem 2.4.18]{Z}. This result will be used later in our analysis.

\begin{pro}
Let $\varphi: \mathbb{R}^n\rightarrow \mathbb{R}$ be a convex and continuous function. Define
$\varphi_+(x):=\max\{\varphi(x), 0\}$ for all $x\in\mathbb{R}^n$. Then $\varphi_+$ is a convex continuous function and
\begin{equation}
  \partial \varphi_+(x)=[0, 1]\partial \varphi(x)
\end{equation}
holds for all $x\in\mathbb{R}^n$ with $\varphi(x)=0$, where $[0, 1]\partial \varphi(x):=\{t\gamma: t\in [0, 1]\ {\rm and}\ \gamma\in\partial \varphi(x)\}$ for any $x\in\mathbb{R}^n$.
\end{pro}

\setcounter{equation}{0}

\section{Main Results}
In this section, we mainly study convex MINLP problem of \eqref{1.1} by dropping the differentiability assumption and aim to establish one outer approximation algorithm for solving such problem.

Let convex MINLP be defined as \eqref{1.1} and set $g:=(g_1,\cdots,g_m)$. For the case that $f,g_i(i=1,\cdots,m)$ in \eqref{1.1} are convex and smooth, it is known from \cite{B1,EMW,FL} that main idea of outer approximation algorithm for convex smooth MINLPs is using linearization of the objective function and the constraints at different points to build a mixed-integer linear program (MILP) relaxation of the problem; that is, given some set $K$ with optimal solutions of several optimization problems, it is possible to build a relaxation of problem (P) in \eqref{1.1}:
\begin{equation}\label{3.1}
\left \{
\begin{array}l
\mathop{\rm minimize}\ \ \ \theta\\
 {\rm subject\ to}\ \ f(x_j, y_j)+\triangledown f(x_j,y_j)^T\begin{pmatrix}x-x_j\\y-y_j\end{pmatrix}\leq \theta,  \\
 \ \ \ \ \ \ \ \ \ \ \ \ \ \ \ g(x_j, y_j)+\triangledown g(x_j,y_j)^T\begin{pmatrix}x-x_j\\y-y_j\end{pmatrix}\leq 0,\\
\ \ \ \ \ \ \ \ \ \ \ \ \ \ \   x\in X, y\in Y\ {\rm discrete \ variable}.
\end{array}
\right. \ \ \forall (x_j,y_j)\in K
\end{equation}

When dealing with problem (P) in \eqref{1.1}, the concept of subgradient is the substitute of the gradient in relaxation of (P). Note that  arbitrary subgradients substituting gradients in \eqref{3.1} is not sufficient to equivalently reformulate problem (P) (see Example 3.1 below). As in \cite{B,EMW}, with the help of KKT conditions, we obtain several special subgradients, which we then use to reformulate problem (P) as one equivalent MILP master program such as \eqref{3.1}.

{\bf{3.1. An overview of the method.}} For the equivalent reformulation of problem (P) in \eqref{1.1} and by using techniques in \eqref{3.1}, we appeal to the concept of projection for expressing problem (P) onto $y$ variables. For any fixed $y\in Y$, we consider the following subproblem $P^y$:
\begin{equation}\label{3.2}
P^y \left \{
\begin{array}l
\mathop{\rm minimize}\limits_{x} \ \ \ f(x, y)\\
 {\rm subject\ to}\ \  g(x, y)\leq 0, \\
\ \ \ \ \ \ \ \ \ \ \ \ \ \ \   x\in X.
\end{array}
\right.
\end{equation}
If there exists some $x\in X$ such that $g(x, y)\leq 0$, the subproblem $P^y$ is said to be feasible; otherwise, $P^y$ is said to be infeasible.

For the validness of KKT conditions, we assume that the following Slater constraint qualification holds:

Assumption (A1) {\it For any $y\in Y$ satisfying that the subproblem $P^y$ is feasible, the following Slater constraint qualification holds:}
$$
g(\hat x, y)<0\ \  {\it for\ some}\ \hat x\in X.  \leqno{\rm (Slater \ CQ)}
$$

Let
\begin{equation}\label{3.3}
    \Sigma:=\{y\in Y: g(x, y)\leq 0\ \ {\rm for\ some}\ x\in X\}
\end{equation}
denote the set of all discrete variables $y$ that produce feasible subproblems. Then the projection of problem (P) onto variable $y$ can be given as follows:
\begin{equation}\label{3.4}
\mathop{\rm minimize}\limits_{y_j\in \Sigma}\left\{
\begin{array}l
\mathop{\rm minimize}\limits_{x} \ \ \ f(x, y_j)\\
 {\rm subject\ to} \ \  g(x, y_j)\leq 0,\\
\ \ \ \ \ \ \ \ \ \ \ \ \ \ \ x\in X.
\end{array}
\right.
\end{equation}

Now let $y_j\in\Sigma$ be fixed. Since $X$ is compact and $f,g_i$ are continuous, it follows that the optimal solution to subproblem $P^{y_j}$ exists. Thus we can suppose that $x_j$ is one optimal solution to $P^{y_j}$. By the assumption (A1) and KKT conditions, there exist $(\lambda_{j,1},\cdots,\lambda_{j,m})\in\mathbb{R}_+^m$ such that
\begin{equation}\label{3.5}
\left\{\begin{array}l
0\in\partial f(\cdot,y_j)(x_j)+\sum\limits_{i\in I(x_j)}\lambda_{j,i}\partial g_i(\cdot,y_j)(x_j)+ N(X, x_j),\\
\lambda_{j,i}g_i(x_j,y_j)=0, \ \  i=1,\cdots,m,\\
\lambda_{j,i}\geq 0, \ \ i=1,\cdots,m,
\end{array}
\right.
\end{equation}
where
\begin{equation}\label{3.6}
I(x_j):=\{i\in\{1,\cdots,m\}: g_i(x_j,y_j)=0\}
\end{equation}
is the active constraint set. This means that we can take $\alpha_j\in\partial f(\cdot,y_j)(x_j)$ and $\xi_{j,i}\in\partial g_i(\cdot,y_j)(x_j) (i=1,\cdots,m)$ such that
\begin{equation}\label{3.7}
  -\alpha_j-\sum\limits_{i\in I(x_j)}\lambda_{j,i}\xi_{j,i}\in N(X, x_j).
\end{equation}
By Proposition 2.3, there exist $\beta_j\in\mathbb{R}^p$ and $\eta_{j,i}\in\mathbb{R}^p (i=1,\cdots,m)$ such that
\begin{equation}\label{3.8}
  (\alpha_j, \beta_j)\in\partial f(x_j,y_j)\ \ {\rm and} \  \ (\xi_{j,i},\eta_{j,i})\in\partial g_i(x_j,y_j), \forall i\in\{1,\cdots,m\}.
\end{equation}
Set $\xi_j:=(\xi_{j,1},\cdots,\xi_{j,m})$ and $\eta_j:=(\eta_{j,1},\cdots,\eta_{j,m})$. We consider the following linear problem:
\begin{equation}\label{3.9}
LP(x_j, y_j) \left \{
\begin{array}l
\mathop{\rm minimize}\limits_{x} \ \ \ f(x_j, y_j)+(\alpha_j, \beta_j)^T\begin{pmatrix}x-x_j\\0\end{pmatrix}\\
 {\rm subject\ to}\ \  g(x_j, y_j)+(\xi_j, \eta_j)^T\begin{pmatrix}x-x_j\\0\end{pmatrix}\leq 0, \\
\ \ \ \ \ \ \ \ \ \ \ \ \ \ \   x\in X.
\end{array}
\right.
\end{equation}

The following theorem establishes the equivalence between subproblem $P^{y_j}$ and linear program $LP(x_j,y_j)$ of \eqref{3.9}. The proof of this theorem will be given in Section 5.

\begin{them}
Let $LP(x_j, y_j)$ be defined as \eqref{3.9}. Then $x_j$ is one optimal solution for $LP(x_j, y_j)$ in \eqref{3.9} and $f(x_j, y_j)$ is the optimal value of $LP(x_j, y_j)$ in \eqref{3.9}.
\end{them}

We denote
\begin{equation}\label{3.13}
T:=\Big\{j: P^{y_j}\ {\rm is \ feasible\ and } \ x_j\ {\rm is\ an \ optimal\ solution\ to}\ P^{y_j}\Big\}.
\end{equation}
Let $j\in T$. By assumption (A1), we can take $(\lambda_{j,1},\cdots,\lambda_{j,m})\in\mathbb{R}^m_+$, $(\alpha_j, \beta_j)\in\partial f(x_j,y_j)$ and $(\xi_{j,i},\eta_{j,i})\in\partial g_i(x_j,y_j) $ $(i=1,\cdots,m)$ such that \eqref{3.7} holds. Applying Proposition 2.3, there exist $\beta_j\in\mathbb{R}^p$ and $\eta_{j,i}\in\mathbb{R}^p$ $(i=1,\cdots,m)$ such that \eqref{3.8} holds. Then we set
$$
\xi_j:=(\xi_{j,1},\cdots,\xi_{j,m}) \ \ {\rm and} \ \ \eta_j:=(\eta_{j,1},\cdots,\eta_{j,m}).
$$
We consider the following MILP:
\begin{equation}\label{3.14}
(M_{\Sigma})\left \{
\begin{array}l
\mathop{\rm minimize}\limits_{x,\,y,\,\theta}\ \ \ \theta\\
 {\rm subject\ to}\ \ f(x_j, y_j)+(\alpha_j, \beta_j)^T\begin{pmatrix}x-x_j\\y-y_j\end{pmatrix}\leq \theta\ \ \forall j\in T,\ \  \\
 \ \ \ \ \ \ \ \ \ \ \ \ \ \ \ g(x_j, y_j)+(\xi_j, \eta_j)^T\begin{pmatrix}x-x_j\\y-y_j\end{pmatrix}\leq 0\ \ \forall j\in T,\\
\ \ \ \ \ \ \ \ \ \ \ \ \ \ \   x\in X, y\in \Sigma\ {\rm discrete \ variable}.
\end{array}
\right.
\end{equation}

By virtue of Theorem 3.1, we obtain the following theorem on the equivalence of problem (P) of \eqref{1.1} and MILP $(M_{\Sigma})$ of \eqref{3.14}.
\begin{them}
Assmue that MINLP problem (P) of \eqref{1.1} satisfies assumption (A1). Then MILP $(M_{\Sigma})$ of \eqref{3.14} are equivalent to problem (P) in the sense that both have the same optimal value and that the optimal solution $(\bar x, \bar y)$ to problem (P) corresponds to the optimal solution $( \bar x, \bar y, \bar\theta)$ to $(M_{\Sigma})$ of \eqref{3.14} with $\bar\theta=f(\bar x, \bar y)$.
\end{them}

For completely reformulating the problem (P), it remains to provide an appropriate representation of constraint $y\in Y\backslash\Sigma$ by supporting hyperplanes. Along the lines in \cite{B1,FL}, we are inspired to study infeasible subproblems so as to eliminate those discrete variables that give rise to infeasibility.

Let $y_l\in Y\backslash\Sigma$. Then subproblem $P^{y_l}$ is infeasible; that is,
$$
\not\exists x\in X \ \ {\rm satisfying} \ g_i(x, y_l)\leq 0\ \ {\rm for \ all} \ i=1,\cdots, m.
$$
Let $J_l$ be one subset of $\{1,\cdots,m\}$ such that there is some $\hat x\in X$ satisfying
\begin{equation}\label{3.24}
  g_i(\hat x, y_l)< 0,\ \ \forall i\in J_l.
\end{equation}
Denote $J_l^{\bot}:=\{1,\cdots,m\}\backslash J_l$ the complement of $J_l$. To detect the infeasibility, we study the following subproblem $F^{y_l}$:
\begin{equation}\label{3.25}
F^{y_l}\left \{
\begin{array}l
\mathop{\rm minimize}\limits_{x} \ \ \ \sum\limits_{i\in J_l^{\bot}}[g_i(x, y_l)]_+\\
 {\rm subject\ to}\ \ g_i(x, y_l)\leq 0\ \ \forall i\in J_l,\\
\ \ \ \ \ \ \ \ \ \ \ \ \ \ \   x\in X,
\end{array}
\right.
\end{equation}
where $[g(x, y_l)]_+:=\max\{g(x, y_l), 0\}$.

Since $X$ is compact and $g_i$ for $i=1,\cdots,m$ are continuous, then the optimal solution to subproblem $F^{y_l}$ exists. Thus we can assume that $x_l$ is one optimal solution to subproblem $F^{y_l}$. For convenience to state the process, we divide the set $J_l^{\bot}$ into three disjoint subsets which are denoted by $J_l^1, J_l^2$ and $J_l^3$. These three subsets are defined as
\begin{equation}
\left\{
\begin{array}l
  J_l^1:=\{i\in J_l^{\bot}: g_i(x_l,y_l)=0\},\\
  J_l^2:=\{i\in J_l^{\bot}: g_i(x_l,y_l)>0\},\\
  J_l^3:=\{i\in J_l^{\bot}: g_i(x_l,y_l)<0\}.
  \end{array}
\right.
\end{equation}
This means that $J_l^{\bot}=J_l^1\cup J_l^2\cup J_l^3$ and by using continuity of $g_i$, one has
$$
\partial[g_i(\cdot, y_l)]_+(x_l)=\partial g_i(\cdot, y_l)(x_l),\ \forall i\in J_l^2\ \ {\rm and} \ \ \partial [g_i(\cdot, y_l)]_+(x_l)=\{0\},\ \forall i\in J_l^3.
$$

By \eqref{3.24} and KKT condition, there exist $\lambda_{l,i}\in\mathbb{R}_+$ for all $i\in J_l$ such that
\begin{equation}\label{3.26}
\left\{
\begin{array}l
0\in\sum\limits_{i\in J_l^{\bot}}\partial[g_i(\cdot,y_l)]_+(x_l)+\sum\limits_{i\in J_l}\lambda_{l,i}\partial g_i(\cdot,y_l)(x_l)+ N(X, x_l),\\
\lambda_{l,i}g_i(x_l,y_l)=0, \ \ \forall i\in J_l,\\
\lambda_{l,i}\geq 0, \ \ \forall i\in J_l.
\end{array}
\right.
\end{equation}
Denote $\lambda_{l,i}\equiv 1$ for all $i\in J_l^2$ and $\lambda_{l,i}\equiv 0$ for all $i\in J_l^3$. Using Proposition 2.4, there exist $\lambda_{l,i}\in [0,1] (\forall i\in J_l^1)$ and $\xi_{l,i}\in\partial g_i(\cdot,y_l)(x_l) (\forall i\in J_l^{\bot}\cup J_l)$ such that
\begin{equation}\label{3.27}
  -\sum_{i\in J_l^{\bot}\cup J_l}\lambda_{l,i}\xi_{l,i}\in N(X,x_l).
\end{equation}
By virtue of Proposition 2.3, there exist $\eta_{l,i}\in\mathbb{R}^p$ such that $(\xi_{l,i}, \eta_{l,i})\in\partial g_i(x_l,y_l)$ for all $i\in J_l^{\bot}\cup J_l$.

Since subproblem $P^{y_l}$ is infeasible, then there exists one optimal solution $x_l$ to subproblem $F^{y_l}$ such that $\sum_{i\in J_l^{\bot}}[g_i(x_l, y_l)]_+>0$, by the continuity of $g_i$ and compactness of $X$. This gives the following theorem on subproblem $F^{y_l}$. The proof is also given in Section 5.

\begin{them}
The discrete variable $y_l\in Y\backslash\Sigma$ is infeasible to the following constraint:
\begin{equation}\label{3.28}
\left \{
\begin{array}l
g_i(x_l, y_l)+(\xi_{l,i}, \eta_{l,i})^T\begin{pmatrix}x-x_l\\y-y_l\end{pmatrix}\leq 0, \ \ \forall i\in J_l^{\bot}\cup J_l,\\
x\in X, y\in Y.
 \end{array}
\right.
\end{equation}
\end{them}

It is necessary to ensure that discrete variables that produce infeasible subproblems are also infeasible in the reformulated master program. We denote
\begin{equation}\label{3.32}
  S:=\big\{l: P^{y_l}\ {\rm is \ infeasible\ and}\ x_l\ {\rm solves}\ F^{y_l} \big\}.
\end{equation}
For any $l\in S$, take $\lambda_{l,i}\geq 0$ and $\xi_{l,i}\in\partial g_i(\cdot,y_l)(x_l) (i=1,\cdots,m)$ such that \eqref{3.27} holds. Take $\eta_{l,i}\in\mathbb{R}^p$ such that $(\xi_{l,i},\eta_{l,i})\in\partial g_i(x_l,y_l)$  for any $i\in\{1,\cdots,m\}$ by Proposition 2.3. We set $\xi_l:=(\xi_{l,1},\cdots,\xi_{l,m})$ and $\eta_l:=(\eta_{l,1},\cdots,\eta_{l,m})$. Then by using Theorem 3.3, we have the following theorem which shows how to eliminate those discrete variables giving rise to infeasible subproblems.

\begin{them}
For any $l\in S$, let $(\xi_l, \eta_l)$ be defined as above. Then the following constraints
\begin{equation}
\left \{
\begin{array}l
g(x_l, y_l)+(\xi_{l}, \eta_{l})^T\begin{pmatrix}x-x_l\\y-y_l\end{pmatrix}\leq 0, \ \ \forall l\in S,\\
x\in X, y\in Y
 \end{array}
\right.
\end{equation}
exclude all discrete variables $y_l\in Y$ for which subproblem $P^{y_l}$ is infeasible.
\end{them}

It is known from Theorem 3.4 that we can add linearization from $F^{y_l}$ when subproblem $P^{y_l}$ is infeasible so as to correctly represent the constraints $y\in\Sigma$ in \eqref{3.3}. This gives rise to the MILP master program (MP) which is equivalent to MINLP problem (P) in \eqref{1.1} and used to reformulate problem (P).

Let $T$ and $S$ be defined as \eqref{3.13} and \eqref{3.32}, respectively. For any $j\in T$, by assumption (A1), we can take $\lambda_{j,i}\geq 0$ $(i=1,\cdots,m)$, $\alpha_j\in\partial f(\cdot,y_j)(x_j)$ and $\xi_{j,i}\in\partial g_i(\cdot,y_j)(x_j) (i=1,\cdots,m)$ such that \eqref{3.7} holds, and by Proposition 2.3, we take $\beta_j\in\mathbb{R}^p$ and $\eta_{j,i}\in\mathbb{R}^p (i=1,\cdots,m)$ such that \eqref{3.8} holds. We set
$$
\xi_j:=(\xi_{j,1},\cdots,\xi_{j,m})\ \ {\rm and} \ \ \eta_j:=(\eta_{j,1},\cdots,\eta_{j,m}).
$$
For any $l\in S$,  we take $\lambda_{l,i}\geq 0$ and $\xi_{l,i}\in\partial g_i(\cdot,y_l)(x_l) (i=1,\cdots,m)$ such that \eqref{3.27} holds and by Proposition 2.3, we take $\eta_{l,i}\in\mathbb{R}^p$ such that $(\xi_{l,i},\eta_{l,i})\in\partial g_i(x_l,y_l)$. Set $$\xi_l:=(\xi_{l,1},\cdots,\xi_{l,m})\ \ {\rm and} \ \ \eta_l:=(\eta_{l,1},\cdots,\eta_{l,m}).
$$
The MILP master problem (MP) is given as follows:
\begin{equation}\label{3.34}
{\rm (MP)}\left \{
\begin{array}l
\mathop{\rm minimize}\limits_{x,\,y,\,\theta} \ \  \ \theta\\
 {\rm subject\ to}\ \ f(x_j, y_j)+(\alpha_j, \beta_j)^T\begin{pmatrix}x-x_j\\y-y_j\end{pmatrix}\leq \theta\ \ \forall j\in T,\ \  \\
 \ \ \ \ \ \ \ \ \ \ \ \ \ \ \ g(x_j, y_j)+(\xi_j, \eta_j)^T\begin{pmatrix}x-x_j\\y-y_j\end{pmatrix}\leq 0\ \ \forall j\in T,\\
 \ \ \ \ \ \ \ \ \ \ \ \ \ \ \ g(x_l, y_l)+(\xi_l, \eta_l)^T\begin{pmatrix}x-x_l\\y-y_l\end{pmatrix}\leq 0\ \ \forall l\in S,\\
\ \ \ \ \ \ \ \ \ \ \ \ \ \ \   x\in X, y\in Y\ {\rm discrete \ variable}.
\end{array}
\right.
\end{equation}

The following theorem, immediate from Theorems 3.3 and 3.4, is one main result in the procedure of reformulating MINLP problem (P) of \eqref{1.1} as the equivalent MILP master program (MP).
\begin{them}
Assume that MINLP problem (P) of \eqref{1.1} satisfies assumptions (A1). Then master program (MP) of \eqref{3.34} is equivalent to problem (P) in the sense that both problems have the same optimal value and that the optimal solution $(\bar x, \bar y)$ to problem (P) corresponds to the optimal solution $(\bar x, \bar y, \bar\theta)$ to (MP) of \eqref{3.34} with $\bar\theta=f(\bar x, \bar y)$.
\end{them}

\noindent{\bf Remark 3.1.} Theorem 3.5 is one extension of main results given in \cite{FL,WA}, and it generalizes the outer approximation method in the sense of equivalently reformulating convex MINLP problem (P) from differentiable case to the non-differentiable one. Further, it is known from Theorem 3.5 that all optimal solutions of problem (P) are optimal solutions to master program (MP). However, the converse is not necessarily true since some optimal solutions of (MP) may be infeasible to problem (P). We refer the reader to \cite[Example 1]{B1} and \cite[Remark 3.1]{WA} for the detail.\\

Theorem 3.5 shows that some subgradients obtained from the KKT conditions enable to reformulate MINLP problem (P) as an equivalent MILP master program by outer approximation method. However, this procedure is not valid if arbitrary subgradients are chosen to replace gradients. The following example demonstrates that the substitution of the gradient by an arbitrary subgradient in the outer approximation method is insufficient for the equivalent reformulation.\\

\noindent{\bf Example 3.1.} We consider the following convex MINLP problem:
\begin{equation}\label{3-21}
\left \{
\begin{array}l
\mathop{\rm minimize}\limits_{x,\,y} \ \ \ f(x, y):=x+y\\
 {\rm subject\ to}\ \  g_1(x, y):=\max\{-x+y+1, x-y+1\}\leq 0,\\
 \ \ \ \ \ \ \ \ \ \ \ \ \ \ \ g_2(x,y):=x-y\leq 0,\\
\ \ \ \ \ \ \ \ \ \ \ \ \ \ \   x\in [0, 2], \ y\in \{1,2,3\}.
\end{array}
\right.
\end{equation}
One can verify that this convex MINLP in \eqref{3-21} is infeasible. However, let the initial point $y_0=1$. Then subproblem $P^{y_0}$ is infeasible and we can consider the following subproblem $F^{y_0}$:
\begin{equation}\label{22}
\left \{
\begin{array}l
\mathop{\rm minimize}\limits_{x} \ \ \ [g_1(x, y_0)]_+=g_1(x, y_0)\\
 {\rm subject\ to}\ \  g_2(x, y_0)\leq 0,\\
\ \ \ \ \ \ \ \ \ \ \ \ \ \ \   x\in [0, 2].
\end{array}
\right.
\end{equation}
It is easy to verify that $x_0=1$ is the optimal solution to subproblem $F^{y_0}$ and
$$
[-1, 1]\times[-1, 1]\subset \partial g_1(x_0, y_0).
$$
Now, if we take $(\xi_{0,1},\eta_{0,1})=(1,1)\in\partial g_1(x_0, y_0)$, $(\xi_{0,2},\eta_{0,2})=\triangledown g_2(x_0, y_0)$ and $(\alpha_0, \beta_0)=\triangledown f(x_0, y_0)$, then the MILP $LP(x_0, y_0)$ is defined as
\begin{equation}\label{23}
\left \{
\begin{array}l
\mathop{\rm minimize}\limits_{x,\,y,\,\theta} \ \  \ \theta\\
 {\rm subject\ to}\ \ x+y\leq \theta,\ \  \\
 \ \ \ \ \ \ \ \ \ \ \ \ \ \ \ x+y-1\leq 0,\\
 \ \ \ \ \ \ \ \ \ \ \ \ \ \ \ x-y\leq 0,\\
\ \ \ \ \ \ \ \ \ \ \ \ \ \ \   x\in [0, 2], \ y\in \{1,2,3\}.
\end{array}
\right.
\end{equation}
The optimal solution to MILP in \eqref{23} is $(x,y,\theta)=(0,1,1)$.
This means that the outer approximation method for this MINLP may generate an infinite loop between points $(x_0,y_0)$ and $(0,1)$. Thus the outer approximation method is invalid for this MINLP problem in \eqref{3-21}. Further, when tracking down why this method is not valid here, the reason noticed is that the KKT conditions at $(x_0,y_0)$ for $(\xi_{0,1},\eta_{0,1})$ does not hold; that is,
$$
\not\exists(\lambda_{0,1},\lambda_{0,2})\in\mathbb{R}^2_+\ \  {\rm satisfying}\ \ \triangledown f(x_0, y_0)+\lambda_{0,1}(\xi_{0,1},\eta_{0,1})+\lambda_{0,2}\triangledown g_2(x_0, y_0)=0.
$$

{\bf 3.2. The algorithm.}
In this subsection, based on the solution of MILP master program (MP) in \eqref{3.34}, we pay main attention to one outer approximation algorithm for finding the optimal solution of problem (P) in \eqref{1.1} along the line in \cite{DG,FL,WA}.


At iteration $k$, the sets $T$ and $S$ in master program (MP) of \eqref{3.34} are substituted by the sets $T^k$ and $S^k$ respectively which are defined as
\begin{equation}\label{24}
\left\{
\begin{array}l
T^k:=\{j\leq k: P^{y_j}\ {\rm is \ feasible \ and} \ x_j \ {\rm solves }\ P^{y_j} \},\\
S^k:=\{l\leq k: P^{y_l}\ {\rm is \ infeasible \ and} \ x_l \ {\rm solves } \ F^{y_l} \}.
\end{array}
\right.
\end{equation}
If $k\in T^k$ then $x_k$ solves $P^{y_k}$ and there exist $(\alpha_k,\beta_k)\in \mathbb{R}^n\times\mathbb{R}^p$ and $(\xi_{k,i},\eta_{k,i})\in\mathbb{R}^n\times\mathbb{R}^p$ for all $i=1,\cdots,m$ such that
\begin{equation}\label{3-25}
\left\{
\begin{array}l
 -\alpha_k-\sum\limits_{i=1}^m\lambda_{k,i}\xi_{k,i}\in N(X, x_k)\ \ {\rm for\ some} \ (\lambda_{k,1}\cdots,\lambda_{k,m})\in\mathbb{R}^m_+,\\
 (\alpha_k,\beta_k)\in\partial f(x_k,y_k),\\
 (\xi_{k,i},\eta_{k,i})\in\partial g_i(x_k,y_k),\ \forall i=1,\cdots,m.
 \end{array}
 \right.
\end{equation}
If $k\in S^k$ then $x_k$ solves $F^{y_k}$ and there exist $(\xi_{k,i},\eta_{k,i})\in\mathbb{R}^n\times\mathbb{R}^p$ for all $i=1,\cdots,m$ such that
\begin{equation}\label{3-26}
\left\{
\begin{array}l
 -\sum\limits_{i=1}^m\lambda_{k,i}\xi_{k,i}\in N(X, x_k)\ \ {\rm for\ some} \ (\lambda_{k,1}\cdots,\lambda_{k,m})\in\mathbb{R}^m_+,\\
 (\xi_{k,i},\eta_{k,i})\in\partial g_i(x_k,y_k),\ \forall i=1,\cdots,m.
 \end{array}
 \right.
\end{equation}
Set
$$
\xi_k:=(\xi_{k,1},\cdots,\xi_{k,m})\ \  {\rm and} \ \ \eta_k:=(\eta_{k,1},\cdots,\eta_{k,m}).
$$
To prevent discrete variable assignment $y_j$ (for any $j\in T^k$) from being the solution to the relaxed master program, it is necessary to define $UBD^k:=\min \{f(x_j,y_j): j\in T^k\}$ and add a constraint $\theta<UBD^k$ to the master program. This gives rise to the following relaxed master program $MP^k$:
\begin{equation}\label{4.1a}
MP^k\left \{
\begin{array}l
\mathop{\rm minimize}\limits_{x,\,y,\,\theta}\ \ \ \theta\\
 {\rm subject\ to}\ \ \theta<UBD^k\\
\ \ \ \ \ \ \ \ \ \ \ \ \ \ \ f(x_j, y_j)+(\alpha_j, \beta_j)^T\begin{pmatrix}x-x_j\\y-y_j\end{pmatrix}\leq \theta\ \ \forall j\in T^k,\ \  \\
 \ \ \ \ \ \ \ \ \ \ \ \ \ \ \ g(x_j, y_j)+(\xi_j, \eta_j)^T\begin{pmatrix}x-x_j\\y-y_j\end{pmatrix}\leq 0\ \ \forall j\in T^k,\\
 \ \ \ \ \ \ \ \ \ \ \ \ \ \ \ g(x_l, y_l)+(\xi_l, \eta_l)^T\begin{pmatrix}x-x_l\\y-y_l\end{pmatrix}\leq 0\ \ \forall l\in S^k,\\
\ \ \ \ \ \ \ \ \ \ \ \ \ \ \   x\in X, y\in Y\ {\rm discrete \ variable}.
\end{array}
\right.
\end{equation}
The new discrete variable assignment $y_{k+1}$ can be obtain by solving $MP^k$ and the whole process is repeated iteratively until the relaxed master program is infeasible.

We are now in a position to state the outer approximation algorithm for solving problem (P) in detail.

\begin{algorithm}
\caption{(Outer Approximation Algorithm)}
\begin{algorithmic}[1]
\STATE Initialization. Given an initial $y_0\in Y$, set $T^0=S^0:=\emptyset$, $UBD^0:=\infty$ and let $k:=1$
 \FOR{$k=1,2, \cdots,$}
 \STATE Solve subproblem $P^{y_k}$
 \IF{$P^{y_k}$ is feasible}
 \STATE Choose one solution $x_k$ of $P^{y_k}$

 Choose $(\alpha_k,\beta_k)$ and $\xi_k:=(\xi_{k,1},\cdots,\xi_{k,m})$, $\eta_k:=(\eta_{k,1},\cdots,\eta_{k,m})$ as in \eqref{3-25}

 Set $T^k:=T^{k-1}\cup\{k\}$, $S^k:=S^{k-1}$ and $UBD^k:=\min\{UBD^{k-1}, f(x_k,y_k)\}$

 \ELSE
 \STATE Solve subproblem $F^{y_k}$ and choose one solution $x_k$ of $F^{y_k}$

 Choose $\xi_k:=(\xi_{k,1},\cdots,\xi_{k,m})$, $\eta_k:=(\eta_{k,1},\cdots,\eta_{k,m})$ as in \eqref{3-26}

 Set $S^k:=S^{k-1}\cup\{k\}$, $T^k:=T^{k-1}$ and $UBD^k:=UBD^{k-1}$
\ENDIF

 \STATE Solve the relaxation $MP^k$ and obtain a new discrete variable $y_{k+1}$

  Set $k:=k+1$ and go back to line 3

\ENDFOR
\end{algorithmic}
\end{algorithm}

Under the assumption of finite cardinality of discrete variable subset $Y$, the following theorem shows Algorithm 1 can detect feasibility or infeasibility of problem (P) in \eqref{1.1} and the procedure in Algorithm 1 terminates after a finite steps. The proof is also given in Section 5.

\begin{them}
Suppose that MINLP problem (P) in \eqref{1.1} satisfies assumption (A1) and the cardinality of $Y$ is finite. Then either problem (P) is infeasible or Algorithm 1 terminates in a finite number steps at an optimal value of problem (P).
\end{them}

\section{Conclusions}

This paper is mainly devoted to the study of one convex MINLP in which objective  and constraint functions are continuous and non-differentiable. With no differentiability assumption, subgradients of objective and constraint functions, the substitute of gradients in convex and smooth MINLP, are chosen from the KKT conditions and used to reformulate the MINLP problem as one equivalent mixed-integer linear program. A counterexample shows that the chosen subgradients, if not satisfying KKT conditions, may be invalid for the MILP reformulation, which demonstrates the necessity of KKT conditions in the equivalent reformulation. By solving a finite sequence of subproblems and relaxed MILP problems, one outer approximation algorithm for this convex MINLP is presented to find the optimal solution of the problem. The finite convergence of the algorithm is also proved. The work of this paper is the extension of references \cite{DG,FL,WA} and also generalizes outer approximation method in the sense of dealing with convex MINLP from differentiable case to the non-differentiable one.

\setcounter{equation}{0}
\section{Appendix: proofs of Theorems 3.1, 3.3 and 3.6}

In this section, we present the proofs of several key results in Section 3.

{\it Proof of Theorem 3.1.} To prove Theorem 3.1, it suffices to show that
\begin{equation}\label{3.10}
  (\alpha_j, \beta_j)^T\begin{pmatrix}x-x_j\\0\end{pmatrix}\geq 0,\ \ \forall x\in X\ {\rm with}\ g(x_j, y_j)+(\xi_j, \eta_j)^T\begin{pmatrix}x-x_j\\0\end{pmatrix}\leq 0.
\end{equation}
Let $x\in X$ be such that
$$
g(x_j, y_j)+(\xi_j, \eta_j)^T\begin{pmatrix}x-x_j\\0\end{pmatrix}\leq 0
$$
and let $I(x_j)$ be defined as \eqref{3.6}. Then
\begin{equation}\label{3.11}
\langle\xi_{j,i}, x-x_j\rangle\leq 0, \forall i\in I(x_j).
\end{equation}
By using \eqref{3.7}, there exists $\gamma\in N(X, x_j)$ such that
\begin{equation}\label{3.12}
  \alpha_j+\sum\limits_{i\in I(x_j)}\lambda_{j,i}\xi_{j,i}+\gamma=0.
\end{equation}
Noting that $X$ is convex, it follows that $x-x_j\in T(X,x_j)$.
This together with \eqref{3.11} and \eqref{3.12} implies that
$$
(\alpha_j, \beta_j)^T\begin{pmatrix}x-x_j\\0\end{pmatrix}=\langle\alpha_j, x-x_j\rangle=-\Big\langle\sum\limits_{i\in I(x_j)}\lambda_{j,i}\xi_{j,i}+\gamma, x-x_j\Big\rangle\geq 0.
$$
Hence \eqref{3.10} holds. The proof is complete. $\Box$

{\it Proof of Theorem 3.3.} Since $X$ is compact and $g$ is continuous, then one has
\begin{equation}\label{3.29}
  \sum_{i\in J_l^{\bot}}[g_i(x_l, y_l)]_+>0.
\end{equation}

Suppose to the contrary that there exists $\hat x\in X$ such that $(\hat x, y_l)$ is feasible to the constraint of \eqref{3.28}. Then
\begin{equation}\label{3.30}
g_i(x_l, y_l)+\langle\xi_{l,i},\hat x-x_l\rangle\leq 0,
\ \ \forall i\in J_l^{\bot}\cup J_l.
\end{equation}
 Noting that $\hat x-x_l\in T(X,x_l)$ by the convexity of $X$, it follows from \eqref{3.27} that there exists $\gamma\in N(X, x_l)$ such that
\begin{equation}\label{3.31}
  \sum_{i\in J_l^1\cup J_l^2\cup J_l^3}\lambda_{l,i}\xi_{l,i}+\gamma=0,
\end{equation}
where $\lambda_{l,i}\equiv 1$ for all $i\in J_l^2$ and $\lambda_{l,i}\equiv 0$ for all $i\in J_l^3$. By multiplying \eqref{3.30} by $\lambda_{l,i}$ for any $i\in J_l^1\cup J_l^2\cup J_l^3$, it follows from \eqref{3.31} that
\begin{eqnarray*}
0&\geq&\sum_{i\in J_l^1\cup J_l^2}\lambda_{l,i}g_i(x_l, y_l)+\Big\langle\sum_{i\in J_l^1\cup J_l^2\cup J_l}\lambda_{l,i}\xi_{l,i}+\gamma, \hat x-x_l\Big\rangle\\
&\geq&\sum_{i\in J_l^2}g_i(x_l, y_l)=\sum_{i\in J^{\bot}_l}[g_i(x_l, y_l)]_+
\end{eqnarray*}
as $\hat x-x_l\in T(X, x_l)$ and $\lambda_{l,i}g_i(x_l,y_l)=0 (\forall i\in J_l)$, which contradicts \eqref{3.29}. The proof is complete. $\Box$

{\it Proof of Theorem 3.6.} (This is similar to the proof for \cite[Theorem 4.1]{WA}. For the sake of completeness, we provide the proof in brief.)

For the proof of Theorem 4.1, we first prove that there is no discrete variable in $Y$  generated\ more  than  once by Algorithm 1. Granting this, it follows from the finite cardinality of $Y$ that the termination of Algorithm 1 holds after a finite number steps.

At iteration $k$, let $(\hat x, \hat y, \hat{\theta})$ be an optimal solution to the relaxed master program $MP^k$. By virtue of Theorem 3.4, one can verify that $\hat y\not=y_l$ for all $l\in S^k$. Suppose to the contrary that $\hat y=y_{j_k}$ for some $j_k\in T^k$. Then $(\hat x, y_{j_k}, \hat{\theta})$ solves the relaxed master program $MP^k$ and
\begin{equation}\label{4.2a}
\left \{
\begin{array}l
\hat{\theta}<UBD^k\leq f(x_{j_k}, y_{j_k}),\\
f(x_{j_k}, y_{j_k})+(\alpha_{j_k}, \beta_{j_k})^T\begin{pmatrix}\hat x-x_{j_k}\\0\end{pmatrix}\leq \hat{\theta},\\
g(x_{j_k}, y_{j_k})+(\xi_{j_k}, \eta_{j_k})^T\begin{pmatrix}\hat x-x_{j_k}\\0\end{pmatrix}\leq 0.
 \end{array}
 \right.
\end{equation}
By using \eqref{3.10} in the proof of Theorem 3.1, one has
$$
(\alpha_{j_k}, \beta_{j_k})^T\begin{pmatrix}\hat x-x_l\\0\end{pmatrix}\geq 0.
$$
This and \eqref{4.2a} imply that $f(x_{j_k}, y_{j_k})\leq \hat{\theta}$, which contradicts $\hat{\theta}<f(x_{j_k}, y_{j_k})$ in \eqref{4.2a}. Hence $\hat y\not=y_j$ for all $j\in T_k$. This means that $\hat y$ is distinct from any $y_j$ for all $j\in T^k\cup S^k$.

Now suppose that the relaxed master program $MP^{k}$ is infeasible for some $k$. Then Algorithm 1 terminate at $k$-th step. Let $\rho$ denote the optimal value of MINLP problem (P). If there is some $j\in T^{k-1}$ such that $f(x_j, y_j)=\rho$, then the conclusion holds. Next, we assume that $f(x_j, y_j)>\rho$ for all $j\in T^{k-1}$. Then $UBD^{k-1}>\rho$ and $k\in T^k$ (otherwise, $k\in S^k$, $UBD^{k}=UBD^{k-1}$ by Algorithm 1 and consequently $MP^{k}$ is feasible, a contradiction).  Thus subproblem $P^{y_k}$ is feasible and $f(x_k, y_k)\geq \rho$.

Suppose to the contrary that $f(x_k, y_k)> \rho$. If $f(x_k, y_k)\leq UBD^{k-1}$, then $\rho<f(x_k, y_k)=UBD^{k}$ and $MP^k$ is feasible, a contradiction. If $f(x_k, y_k)> UBD^{k-1}$ then $\rho<UBD^{k-1}=UBD^{k}$ and thus $MP^{k}$ is feasible, a contradiction. This means $f(x_k, y_k)=\rho$. The proof is complete. $\Box$\\

\noindent{\bf Acknowledgment.} We are grateful to the referee for careful reading this paper and valuable comments which help us to improve the original version. This research was supported by the National Natural Science Foundations of P. R. China (Grant No. 11401518 and No. 11261067) and IRTSTYN, and by the Claude Leon Foundation of South Africa.

\end{document}